\documentclass[12pt]{article}
\pdfoutput=1

\usepackage{amsmath}
\usepackage{amssymb}
\usepackage{amsthm}
\usepackage[msc-links]{amsrefs} 
\usepackage{enumerate}

\newcommand{\N}{\mathbb N}

\newcommand{\R}{\mathbb R}
\newcommand{\Z}{\mathbb Z}

\DeclareMathOperator{\rk}{rk}

\newcommand{\mcg}{\mathrm{Mod}}
\newcommand{\spinc}{\(\mathrm{Spin}^c\,\)}

\newcommand{\hfh}{\widehat{HF}}
\newcommand{\hfk}{\widehat{HFK}}

\theoremstyle{plain}
\newtheorem{thm}{Theorem}
\newtheorem{lem}[thm]{Lemma}
\newtheorem{prop}[thm]{Proposition}
\newtheorem{cor}[thm]{Corollary}

\theoremstyle{definition}
\newtheorem{defn}[thm]{Definition}

\theoremstyle{remark}

\theoremstyle{plain}
\newtheorem*{thm_main}{Theorem \ref{thm_main}}
\newtheorem*{prop_pA}{Proposition \ref{prop_pA}}
\newtheorem*{thm_totalrk}{Theorem \ref{thm_totalrk}}

\title{3-Manifold Mutations Detected by {H}eegaard {F}loer Homology}
\author{Corrin Clarkson\thanks{The author was partially supported by NSF grant number DMS-0739392.}}

\begin{document}
\maketitle

\begin{abstract}
	Given a self-diffeomorphism $\varphi$ of a closed, orientable surface $S$ with genus greater than one and an embedding $f$ of $S$ into a three-manifold $M$, we construct a mutant manifold by cutting $M$ along $f(S)$ and regluing by $f\varphi f^{-1}$. We will consider whether there exist nontrivial gluings such that for any embedding, the manifold $M$ and its mutant have isomorphic Heegaard Floer homology. In particular, we will demonstrate that if $\varphi$ is not isotopic to the identity map, then there exists an embedding of $S$ into a three-manifold $M$ such that the rank of the non-torsion summands of $\hfh$ of $M$ differs from that of its mutant.	We will also show that if the gluing map is isotopic to neither the identity nor the genus-two hyperelliptic involution, then there exists an embedding of $S$ into a three-manifold $M$ such that the total rank of $\hfh$ of $M$ differs from that of its mutant.
\end{abstract}

\section{Introduction}

In 2004, Ozsv\' ath and Szab\' o introduced Heegaard Floer homology, a topological invariant that assigns a collection of abelian groups to each closed, oriented three-manifold equipped with a \spinc -structure \citelist{\cite{Ozsvath2004}}. Given a topological invariant, it is natural to ask which topological operations it detects. In this paper, we will consider whether or not Heegaard Floer homology detects \emph{mutation}, the operation of cutting a three-manifold along an embedded surface and regluing by a surface diffeomorphism. In particular, we will show that the version of Heegaard Floer homology denoted by $\hfh$ can detect mutation by any nontrivial diffeomorphisms of a closed, orientable surface of genus greater than one.

In order to make this statement more precise, we introduce the following terminology and notation. Let $g\ge 2$ be a natural number and let $S_g$ be a genus-$g$ smooth, orientable, closed surface. By a \emph{manifold-surface pair}, we will mean a pair $(M, f)$ where $M$ is a closed, connected, smooth 3-manifold and $f\colon S_g \rightarrow M$ is a smooth embedding of $S_g$ into $M$ such that $f(S_g)$ separates $M$. To a diffeomorphism $\varphi\colon S_g\to S_g$ and a manifold-surface pair $(M, f)$, we associate the \emph{mutant manifold} $M_f^\varphi$ that results from cutting $M$ along $f(S_g)$ and regluing by $f\varphi f^{-1}$. 

\begin{thm}\label{thm_main} 
	Let $\varphi$ be a self-diffeomorphism of $S_g$ that is not isotopic to the identity map. Then, there exists a manifold surface pair $(M,f)$ such that 
	$$\mathrm{rk} \bigoplus_{c_1(\mathfrak{s})\ne 0} \hfh(M,\mathfrak{s}) \ne \mathrm{rk} \bigoplus_{c_1(\mathfrak{s})\ne 0} \hfh(M_f^\varphi,\mathfrak{s}).$$
	Here, $c_1(\mathfrak{s})$ is the first Chern class of the \spinc-structure $\mathfrak{s}$.
\end{thm}

Because the manifolds in our manifold-surface pairs are orientable rather than oriented, there is some ambiguity when assigning \spinc structures. However, the rank of $\hfh$ is preserved by both reversing the manifold's orientation and conjugating the \spinc -structure \citelist{\cite{Ozsvath2004b}*{Thm. 2.4 \& Prop. 2.5}}. Thus, the ambiguity is resolved by the fact that we are only interested in the total rank of the non-torsion, \emph{i.e.} $c_1(\mathfrak{s}) \ne 0$ summands.

Our proof of this result begins with a reformulation of the theorem statement. In Section \ref{sec_reform}, we use Ivanov and Long's results about subgroups of mapping class groups to show that Theorem \ref{thm_main} is equivalent to the statement that a particular subgroup of the mapping class group $\mcg(S_g)$ contains neither the genus-2 hyperelliptic involution nor any pseudo-Anosov elements. In Section \ref{sec_genus-2}, we show that the genus-2 hyperelliptic involution is not an element of this subgroup by giving an example of a mutation by this map that changes the rank of the non-torsion summands of $\hfh$. In Section \ref{sec_pA}, we use the fact that $\hfh$ detects the Thurston semi-norm on homology to establish the existence of mutations by pseudo-Anosov maps that change the rank of the non-torsion summands of $\hfh$ \citelist{\cite{Ozsvath2004a} \cite{Ni2009}\cite{Hedden2010}}. We conclude the proof of Theorem \ref{thm_main} in Section \ref{sec_proof}.

Then in Section \ref{sec_totalrk}, we use similar techniques to show that the total rank of $\hfh$ can detect mutations by non-central mapping classes:

\begin{thm}\label{thm_totalrk} 
	Let $[\varphi]\in\mcg(S_g)$ be a mapping class that is isomorphic to neither the identity nor the genus-2 hyperelliptic involution. Then there exists a manifold surface pair $(M,f)$ such that
	$$\mathrm{rk} \ \hfh(M) \ne \mathrm{rk} \ \hfh(M_f^\varphi).$$
\end{thm}

The effect of mutating by the genus-2 hyperelliptic involution has been considered for invariants related to $\hfh$. 
In particular, Ozsv\' ath and Szab\'o showed that the Heegaard Floer knot invariant $\hfk$ can detect mutations of this form \citelist{\cite{Ozsvath2004a}*{Thm. 1.2}}. Conversely, there is computational evidence that the total rank of $\hfk$ is preserved by mutation by the genus-2 hyperelliptic involution \cite{Moore2012}. Finally, Ruberman showed that the instanton Floer homology with $\Z / 2\Z$ coefficients of an oriented homology 3-sphere is preserved by mutations of this form \citelist{\cite{Ruberman1999}*{Thm. 1}}\footnote{In private communication, Ruberman indicated that there is an issue with the signs in this paper due to a particular moduli space not being orientable. However, this is not relevant when one considers $\Z / 2\Z$ coefficients.}. 

The results of this paper also fit into the growing body of work on group actions on triangulated categories. See Section \ref{sec_implications} for a more detailed discussion.

\subsection{Acknowledgments}
I am grateful to Robert Lipshitz for suggesting this problem and for many useful discussions and ideas. I am also grateful to Jason Behrstock, Ian Biringer, Nathan Dunfield, Julia Elisenda Grigsby, Adam Levine, Dan Margalit, Walter Neumann and Dylan Thurston for helpful conversations.

\section{Theorem reformulation} \label{sec_reform}

Let $\mcg(S_g)$ be the mapping class group of $S_g$. In this section, we will reformulate Theorem \ref{thm_main} as a statement about the triviality of a normal subgroup of $\mcg(S_g)$.

\begin{defn} 
	A mapping class $[\varphi]\in\mcg(S_g)$ is \emph{$\hfh$-invisible} if for all manifold-surface pairs $(M,f)$ we have that \[ \rk{\bigoplus_{c_1(\mathfrak{s}) \ne 0}\hfh(M,\mathfrak{s})} = \rk{\bigoplus_{c_1(\mathfrak{s}) \ne 0}\hfh(M_f^\varphi,\mathfrak{s})}.\] 
\end{defn}

The $\hfh$-invisible mapping classes are well defined, because mutating by isotopic diffeomorphisms results in  diffeomorphic mutant manifolds. Moreover, they form a normal subgroup.

\begin{prop} \label{prop_normal} 
	The $\hfh$-invisible mapping classes form a normal subgroup of $\mcg(S_g)$. 

	\begin{proof}
	The mapping class of the identity map is $\hfh$-invisible, because mutating by any of its representatives preserves the diffeomorphism class of the manifold. We will show that mutations by products, inverses and conjugates of $\hfh$-invisible mapping classes preserve the rank of the non-torsion summands of $\hfh$.

	Let $(M, f)$ be a manifold-surface pair and let $M_1$ and $M_2$ be the closures of the two connected components of $M\setminus f(S_g)$. Also let $\alpha$ and $\beta$ be arbitrary self-diffeomorphisms of $S_g$. The mutant manifold $M_f^\alpha$ can be made into a manifold-surface pair by composing the embedding $f|_{M_1}\colon  S_g \rightarrow M_1$ with the inclusion of $M_1$ into $M_f^\alpha$. Let $(N, h)$ denote this pair. Mutating $(N, h)$ by $\beta$ gives the mutant $N^\beta_h$ which is constructed by using $(f\alpha) \beta f^{-1}$ to glue $M_1$ to $M_2$. Thus, $N^\beta_h$ is diffeomorphic to $M^{\alpha\beta}_f$ by construction, and we can view mutation by a composite map as a sequence of mutations. 

Let $[\varphi]$ and $[\tau]$ be an $\hfh$-invisible mapping classes. It follows that mutating by either $\varphi$ or $\tau$ preserves the rank of the non-torsion summands of $\hfh$. Thus, if we view mutating $(M,f)$ by the composition $\varphi\tau$ as a mutation by $\varphi$ followed by a mutation by $\tau$, we find that 
$$\rk{\bigoplus_{c_1(\mathfrak{s}) \ne 0}\hfh(M, \mathfrak{s})} = \rk{\bigoplus_{c_1(\mathfrak{s}) \ne 0}\hfh(M^\varphi_f, \mathfrak{s})}=\rk{\bigoplus_{c_1(\mathfrak{s}) \ne 0}\hfh(M^{\varphi\tau}_f, \mathfrak{s})}.$$
Therefore, the product $[\varphi][\tau] = [\varphi\tau]$ is $\hfh$-invisible.

	Mutating $(M, f)$ by the composite map $\varphi^{-1}\varphi$ does not change it's diffeomorphism class. Furthermore, if we view this mutation sequentially, the second mutation preserves the rank of the non-torsion summands of $\hfh$. Thus, we have that
$$\rk{\bigoplus_{c_1(\mathfrak{s}) \ne 0}\hfh(M, \mathfrak{s})}=   \rk{\bigoplus_{c_1(\mathfrak{s}) \ne 0}\hfh(M_f^{\varphi^{-1}\varphi}, \mathfrak{s})} = \rk{\bigoplus_{c_1(\mathfrak{s}) \ne 0}\hfh(M^{\varphi^{-1}}_f, \mathfrak{s})}.$$
Therefore, the inverse mapping class $[\varphi]^{-1} = [\varphi^{-1}]$ is $\hfh$-invisible.

	Let $[\psi]\in \mcg(S_g)$ be and arbitrary mapping class. Composing $f$ with $\psi$ gives a new embedding $f \psi\colon S_g \rightarrow M$. Mutating the manifold-surface pair $(M, f   \psi)$ by $\varphi$ gives the mutant manifold $M_{f\psi}^\varphi$. This mutant is constructed by using $(f   \psi)   \varphi   (f \psi)^{-1}$ to glue $M_1$ to $M_2$. Similarly, the mutant $M_f^{\psi\varphi\psi^{-1}}$ is constructed by using $f ( \psi   \varphi  \psi^{-1} ) f^{-1}$ to glue $M_1$ to $M_2$ and is thus diffeomorphic to $M_{f\psi}^\varphi$. Furthermore, the rank of the non-torsion summands of $\hfh(M_{f\psi}^\varphi)$ is the same as that of $M$, because $[\varphi]$ is $\hfh$-invisible. Thus, we have that

	$$\rk{\bigoplus_{c_1(\mathfrak{s}) \ne 0}\hfh(M, \mathfrak{s})} =  \rk{\bigoplus_{c_1(\mathfrak{s}) \ne 0}\hfh(M_{f\psi}^\varphi, \mathfrak{s})} = \rk{\bigoplus_{c_1(\mathfrak{s}) \ne 0}\hfh(M_f^{\psi\varphi\psi^{-1}}, \mathfrak{s})}.$$
	Therefore, the conjugate $[\psi][\varphi][\psi]^{-1}=[\psi\varphi\psi^{-1}]$ is $\hfh$-invisible. It follows that the $\hfh$-invisible mapping classes form a normal subgroup of $\mcg(S_g)$.
	\end{proof}
\end{prop} 

Theorem \ref{thm_main} is equivalent to the statement that the normal subgroup of $\hfh$-invisible mapping classes is trivial. Reformulating the theorem statement in this way allows us to leverage the group structure of $\mcg(S_g)$. We begin by recalling a few definitions from the theory of mapping class groups.

The \emph{Torelli group} is the normal subgroup consisting of those mapping classes whose representatives induce the identity map on homology and is denoted by $\mathcal{I}(S_g)$. If $g = 2$, then $\mcg(S_g)$ has a unique order two element that acts by $-id$ on $H_1(S_2;\Z)$ \citelist{\cite{Farb2011}*{\S 7.4}}. This element is called the \emph{genus-2 hyperelliptic involution}. A subgroup $G \le \mcg(S_g)$ is called \emph{irreducible} if for any simple closed curve $C$ on $S_g$ there exists and element $[\varphi] \in G$ such that $\varphi(C)$ is not isotopic to $C$.

We are now ready to state and prove the following proposition:

\begin{prop} \label{prop_MCGtheory} 
	If a normal subgroup $G \vartriangleleft \mcg (S_g)$ contains no pseudo-Anosov elements of the Torelli group, then it is either the trivial subgroup or the order two subgroup generated by the genus-2 hyperelliptic involution.
	
	\begin{proof}
		Let $G \vartriangleleft \mcg (S_g)$ be a normal subgroup of the mapping class group that contains no pseudo-Anosov elements of the Torelli group. Also let $H = G \cap \mathcal{I}(S_G)$ be the intersection of $G$ with the Torelli group. Thus, $H$ is a normal subgroup that contains no pseudo-Anosov elements. 
		
		It follows from a theorem of Ivanov that $H$ is either finite or reducible \citelist{\cite{Ivanov1992}*{Thm. 1}}. Furthermore, the Torelli group is torsion free and thus $H$ must be either trivial or infinite and reducible \citelist{\cite{Ivanov1992}*{Cor. 1.5}}. However, Ivanov also showed that there are no infinite, reducible, normal subgroups of $\mcg(S_g)$ \citelist{\cite{Ivanov1992}*{Cor. 7.13}}. Therefore, $H$ must be trivial.
	
	Long showed that if the intersection of two normal subgroups of $\mcg(S_g)$ is trivial, then one of those groups must either be central or trivial \citelist{\cite{Long1986}*{Lem. 2.1}}. The Torelli group is neither central nor trivial, so we must conclude that $G$ is either central or trivial. If $g \ge 3$, then the center of $\mcg(S_g)$ is trivial \citelist{\cite{Farb2011}*{Thm. 3.10}} and thus $G$ must also be trivial. In the genus-2 case, things are only slightly more complicated. The center of $\mcg(S_2)$ is the order two subgroup generated by the hyperelliptic involution \citelist{\cite{Farb2011}*{\S 3.4}}. Therefore, $G$ is either trivial or the order two subgroup generated by the genus-2 hyperelliptic involution.
	\end{proof}
\end{prop}

By combining Propositions \ref{prop_normal} and  \ref{prop_MCGtheory}, we see that Theorem \ref{thm_main} is equivalent to the statement that neither the genus-2 hyperelliptic involution nor any pseudo-Anosov elements of the Torelli group are $\hfh$-invisible. In the next two sections, we will consider mutations by these two types of mapping classes.

\section{Genus-two hyperelliptic involution} \label{sec_genus-2}

In this section, we will show that mutating by the genus-2 hyperelliptic involution can change the rank of the non-torsion summands of $\hfh$. To accomplish this, we will consider the semi-norm on $H_2(M;\R)$ defined by Thurston in \citelist{\cite{Thurston1986}}. This is a useful invariant to consider, because it is detected by $\hfh$ and is much easier to compute \citelist{\cite{Ozsvath2004a} \cite{Ni2009}\cite{Hedden2010}}. 

\begin{prop} \label{prop_genus-2}
	The genus-2 hyperelliptic involution is not $\hfh$-invisible.

	\begin{proof}
		We consider the pair of mutant knots that form the basis of Moore and Starkston's examples of mutations by the genus-2 hyperelliptic involution \cite{Moore2012}. Let $K$ and $K^\tau$ be the knots denoted respectively by $14^n_{22185}$ and $14^n_{22589}$ in Knotscape notation \citelist{\cite{Moore2012}*{Fig. 2}}. These two knots are related by a mutation of $S^3$ by the genus-2 hyperelliptic involution \citelist{\cite{Moore2012}*{Fig. 3}}. From the computations of $\hfk$ in Table 1 of \cite{Moore2012}, we see that $K$ has genus two and $K^\tau$ has genus one. 
		
		Now, let $M$ and $M^\tau$ be the results of 0-surgery on $K$ and $K^\tau$ respectively. Because the mutation of $S^3$ that transforms $K$ into $K^\tau$ involves a surface that is disjoint from the knot, there is a corresponding surface in $M$. Moreover, mutating $M$ along that corresponding surface by the genus-2 hyperelliptic involution will result in an manifold diffeomorphic to $M^\tau$.
		
		A Mayer-Vietoris argument shows that both $H_2(M;\R)$ and $H_2(M^\tau;\R)$ are isomorphic to $\R$. Furthermore, it follows from the work of Gabai that the genera of the knots $K$ and $K^\tau$ determine the Thurston semi-norm on these homology groups  \citelist{\cite{Gabai1987a}*{Cor. 8.3}}. In particular, the semi-norm is constantly zero on $H_2(M^\tau;\R)$ and nonzero on $H_2(M;\R)$. This implies that  $\hfh(M^\tau)$ is supported entirely in the \spinc -sturcture whose first Chern class is zero and $\hfh(M)$ is nontrivial in at least one \spinc structure with nonzero first Chern class \citelist{\cite{Hedden2010}*{Thm. 2.2}}.
	\end{proof}
\end{prop}

\section{Pseudo-Anosov gluings} \label{sec_pA}

In this section, we examine mutations by pseudo-Anosov elements of the Torelli group. In particular, we will show that mutations of this form can change the Thurston semi-norm on homology.

\begin{prop}\label{prop_pA} 
	Let $[\varphi]\in\mathcal{I}(S_g)$ be a pseudo-Anosov element of the Torelli group. Then, there exists a natural number $N$ and a manifold surface pair $(M,f)$ such that $M = S^1 \times S^2$ and the the mutant manifold $M_f^{\varphi^N}$  has a homology class with nonzero Thurston semi-norm.
\end{prop}

In order to determine the effect of mutation on the Thurston semi-norm, we must first establish a relationship between the homology of a three-manifold and that of its mutants. In the case of mutation by elements of the Torelli group, this is achieved by the following lemma.

\begin{lem}\label{lem_H-invisible}
	If $[\psi] \in \mathcal{I}(S_g)$ is an element of the Torelli group and $(M,f)$ is a manifold surface pair, then $M$ and its mutant $M_f^\psi$ have isomorphic homology groups
	$$H_i\left(M\right) \cong H_i(M_f^\psi) \mbox{ for all } i.$$
	
	\begin{proof}
		Because $M$ and its mutant $M_f^\psi$ are closed three-manifolds, it suffices to show that the first homology groups are isomorphic. In order to do this, we decompose $M$ into two open sets that overlap in a tubular neighborhood of the separating surface $f(S_g)$. A comparison of the Mayer-Vietoris sequence coming from this decomposition to that coming from a similar decomposition of the mutant $M_f^\psi$ shows that the first homology groups are indeed isomorphic.
	\end{proof}
\end{lem}

Our inquiry will focus on mutating $S^1 \times S^2$ along Heegaard surfaces. We proceed by considering the relationship between the complexity of the Heegaard splittings of a three-manifold and the minimal genera of its homology classes.

\subsection{Homology and Hempel Distance}

A genus-$g$ \emph{Heegaard splitting} is a decomposition of a three-manifold into two genus-$g$ handlebodies glued together along their boundaries. Such a splitting is determined by two handlebodies with parameterized boundaries. A handlebody with parameterized boundary is in turn determined by the curves on the boundary that bound disks in the handlebody. 

\begin{defn} 
	For a genus-$g$ handlebody $X$ with boundary parameterized by a map to $S_g$, let $\mathcal{V}_X$ be the set of isotopy classes of essential simple closed curves in $S_g$ whose preimages bound disks in $X$. We will refer to the elements of $\mathcal{V}_X$ as \emph{compression curves} of $X$. 
\end{defn}

Given two genus-$g$ handlebodies $X$ and $Y$ with boundaries parameterized respectively by maps $a$ and $b$ to $S_g$, we can construct a 3-manifold $M$ by using $b^{-1}a\colon  \partial X \to \partial Y$ to glue $X$ to $Y$. We will write $(S_g, \mathcal{V}_X, \mathcal{V}_Y)$ for the corresponding Heegaard splitting of $M$. 

The compression curves of a genus-$g$ handlebody can be viewed as points in the curve complex, $C(S_g)$ \citelist{\cite{Harvey1981}}. The \emph{curve complex} is a simplicial complex with 0-simplicies corresponding to isotopy classes of essential closed curves and n-simplicies corresponding to $(n-1)$-tuples of isotopy classes that can be realized disjointly. There is a natural distance function $d$ on the 0-simplicies of the curve complex given by viewing the 1-skeleton as a graph with edge length one. Applying this distance function to the sets of compression curves in a Heegaard splitting can provide information about the minimal genera of homology classes.

\begin{lem}\label{lem_H-distance} 
	If $(S_g, \mathcal{V}_X, \mathcal{V}_Y)$ is a Heegaard splitting of a manifold $M$ and the distance $d(\mathcal{V}_X, \mathcal{V}_Y)$ is greater than two, then $M$ is irreducible and has no essential tori.
	\begin{proof}
		Haken showed that if $M$ were reducible, then $\mathcal{V}_X$ and $\mathcal{V}_Y$ would have a point in common and thus $d(\mathcal{V}_X, \mathcal{V}_Y)$ would be zero \citelist{\cite{Haken1968}*{pg. 84}}.
Furthermore, Hempel demonstrated that if $M$ had an essential torus, then $d(\mathcal{V}_X, \mathcal{V}_Y)$ would be $\le 2$ \citelist{\cite{Hempel2001}*{Cor. 3.7}}. Thus, $d(\mathcal{V}_X, \mathcal{V}_Y) > 2$ implies that $M$ is irreducible and has no essential tori.
	\end{proof}
\end{lem}

The distance between the sets of compression curves in a Heegaard splitting is called the \emph{Hempel distance} of that splitting. Combining this language with the definition of Thurston's semi-norm gives the following corollary to Lemma \ref{lem_H-distance}.

\begin{cor}\label{cor_norm}
	If a three-manifold $M$ has a Heegaard splitting with Hempel distance greater than two, then the Thurston semi-norm is in fact a norm on $H_2(M;\R)$.
\end{cor}

Now that we have established a relationship between Hempel distance and the Thurston norm, we turn our attention to the effect of mutating by a pseudo-Anosov map on the Hempel distance of a Heegaard splitting.

\subsection{Effects of Pseudo-Anosov mutations}

Each pseudo-Anosov map $\varphi\colon S_g\to S_g$ has two associated projective measured laminations on $S_g$ called its \emph{stable} and \emph{unstable laminations} \citelist{\cite{Casson1988}*{Thm. 5.5}}. Furthermore, a set of compression curves can be viewed as a subset of $PML(S_g)$, the space of projective measured laminations on $S_g$ by simply applying the counting measure to each curve \citelist{\cite{Hamenstadt2007}*{\S 2}}. We will use $\overline{\mathcal{V}_H}$ to denote the closure of $\mathcal{V}_H$ in $PML(S_g)$. Hempel showed that repeatedly twisting by a pseudo-Anosov map can increase the Hempel distance of a Heegaard splitting: 

\begin{thm}[Hempel \citelist{\cite{Hempel2001}*{p. 640}} See also \citelist{\cite{Abrams2005}*{\S 2}}]\label{thm_Hempel}
	
	Let $X$ and $Y$ be genus-$g$ handlebodies with their boundaries parametrized by a maps to $S_g$ and let $\varphi\colon  S_g \to S_g$ be a pseudo-Anosov map with stable lamination $s$ and unstable lamination $u$. If $s$ and $u$ are not in $\overline{\mathcal{V}_X} \cup \overline{\mathcal{V}_Y}$, then the distance between $\mathcal{V}_X$ and $\varphi^n(\mathcal{V}_Y)$ tends to infinity,
	\[\lim_{n\rightarrow \infty} d\left(\mathcal{V}_X,\varphi^n(\mathcal{V}_Y)\right) = \infty.\]
\end{thm}

It is worth noting that $(S^g, \mathcal{V}_X,\varphi^n(\mathcal{V}_Y))$ is a the Heegaard splitting of the mutant manifold that results from mutating $X \cup Y$ by $\varphi^n$ along the Heegaard surface $\partial X$. We would like to use Hempel's theorem to make statements about mutations of $S^1\times S^2$ by pseudo-Anosov maps. However, we must first verify that $S^1\times S^2$ admits Heegaard splittings of the appropriate form. 

\begin{lem}\label{lem_hbody}
	Let $\varphi\colon S_g\to S_g$ is a pseudo-Anosov map with stable lamination $s$ and unstable lamination $u$. Then there exists a genus-$g$ Heegaard splitting $(S_g, \mathcal{V}_X, \mathcal{V}_Y)$ of $S^1\times S^2$ such that $s$ and $u$ are not in $\overline{\mathcal{V}_X} \cup \overline{\mathcal{V}_Y}$.
	
	\begin{proof}	
	For an arbitrary handlebody $X$, the stable lamination $s$ is in $\overline{\mathcal{V}_X}$ if and only if the unstable lamination $u$ is also in $\overline{\mathcal{V}_X}$ \citelist{\cite{Biringer2010}*{Thm. 1.1}}. Thus, it is enough to find a Heegaard splitting of $S^1\times S^2$ such that $s$ is not in the closure of either set of compression curves.
	
	Let $(S_g, \mathcal{V}_X, \mathcal{V}_Y)$ be a genus-$g$ Heegaard Splitting of $S^1\times S^2$. The union $\overline{\mathcal{V}_X} \cup \overline{\mathcal{V}_Y}$ is nowhere dense in $PML(S_g)$\citelist{\cite{Masur1986}*{Thm. 1.2}}. Furthermore, The stable laminations of psuedo-Anosov elements of $\mcg(S_g)$ form a dense subset of $PML(S_g)$ \citelist{\cite{Fathi2012}*{Thm. 6.19}}. Thus, there exists a pseudo-Anosov map $\psi:S_g\to S_g$ with stable lamination $t$ such that $t$ is not in $\overline{\mathcal{V}_X} \cup \overline{\mathcal{V}_Y}$ and $t$ is not equal to $s$ or $u$.
	
	We will now show that translating the set $\mathcal{V}_X$ by a high power of $\psi$ will move it away from $s$. By Theorem \ref{thm_Hempel}, we have that for any $n\in \N$ the distance  $d(\psi^n(\mathcal{V}_X), \psi^{n+m}(\mathcal{V}_X))$ goes to infinity as $m$ grows. Thus, it is enough to show that if $s$ is a limit point of both $\psi^n(\mathcal{V}_X)$ and $\psi^{n+m}(\mathcal{V}_X)$ in $PML(S_g)$, then these sets must be close together in the curve complex.

	Suppose $s$ is an element of both $\overline{\psi^n(\mathcal{V}_X)}$ and $\overline{\psi^{n+m}(\mathcal{V}_X)}$. Let $(a_i)$ and $(b_i)$ be sequences of points in $\psi^n(\mathcal{V}_X)$ and $\psi^{n+m}(\mathcal{V}_X)$ respectively that converge to $s$ in $PML(S_g)$. It follows from work of Klarreich that the sequences $(a_i)$ and $(b_i)$ converge to the same point in the Gromov boundary of the curve complex $C(S_g)$ \citelist{\cite{Klarreich1999}*{}}. See also \citelist{ \cite{Abrams2005}*{Thm. 8.4}} and \citelist{\cite{Hamenstadt2006}*{Thm. 1}}. 
This in turn implies that the Hempel distance between $\psi^{n}(\mathcal{V}_X)$ and $\psi^{n+m}(\mathcal{V}_X)$ is bounded above by a constant $K$ which depends only on the genus $g$ \citelist{\cite{Abrams2005}*{Lem. 9.2}}. 

	Therefore, there exists an $M\in \N$ such that $s$ is not in $\overline{\psi^n(\mathcal{V}_X)}$ for all $n >M$. Similarly, translating $\mathcal{V}_Y$ by a high power of $\psi$ will move it away from $s$. Thus, there exists an $N$ such that $s$ is not in $\overline{\psi^N(\mathcal{V}_X)} \cup \overline{\psi^N(\mathcal{V}_Y)}$. By construction,  $(S_g, \psi^N(\mathcal{V}_X), \psi^N(\mathcal{V}_Y))$ is a Heegaard splitting for $S^1\times S^2$.
	\end{proof}
\end{lem}


\subsection{Proof of Proposition \ref{prop_pA}}

\begin{prop_pA} 
	Let $[\varphi]\in\mathcal{I}(S_g)$ be a pseudo-Anosov element of the Torelli group. Then, there exists a natural number $N$ and a manifold surface pair $(M,f)$ such that $M = S^1 \times S^2$ and the the mutant manifold $M_f^{\varphi^N}$ has a homology class with nonzero Thurston semi-norm.
\end{prop_pA}

\begin{proof}
	Let $s, u \in PML(S_g)$ be respectively the stable and unstable laminations of $\varphi$. Also, let $(S_g, \mathcal{V}_X, \mathcal{V}_Y)$ be a genus-$g$ Heegaard splitting of $S^1\times S^2$ such that $s$ and $u$ are not in $\overline{\mathcal{V}_X}\cup \overline{\mathcal{V}_Y}$. The existence of such a splitting is guaranteed by Lemma \ref{lem_hbody}. Finally, let $(M,f)$ be the manifold surface pair where $M=S^1\times S^2$ and $f$ is the embedding of $S_g$ as the Heegaard surface $\partial X$ from the splitting $(S_g, \mathcal{V}_X, \mathcal{V}_Y)$.
	
	By Theorem \ref{thm_Hempel}, we have that 
	\[\lim_{n\rightarrow \infty} d\left(\mathcal{V}_X,\varphi^n(\mathcal{V}_Y)\right) = \infty.\] 
	Thus, there exists a natural number $N$ such that $d\left(\mathcal{V}_X,\varphi^N(\mathcal{V}_Y)\right) > 2$. Furthermore, $\left(S_g, \mathcal{V}_X,\varphi^N(\mathcal{V}_Y)\right)$ is a Heegaard splitting for the mutant $M_f^{\varphi^N}$. This implies that $M_f^{\varphi^N}$ is irreducible and has no essential tori (Lem. \ref{lem_H-distance}). 
	
	A simple calculations shows that the $H_2(M;\Z) = H_2(S^1\times S^2;\Z) \cong \Z$. It follows that $H_2(M^{\varphi^N}_f;\Z) \cong \Z$, because $[\varphi]$ is in the Torelli group (Lem. \ref{lem_H-invisible}). Let $\omega$ be a nonzero element of $H_2(M^{\varphi^N}_f;\Z) \cong \Z$ and let $F\subseteq M^{\varphi^N}_f$ be a surface that represents $\omega$. Because $M^{\varphi^N}_f$ is irreducible and has no essential tori, the genus of $F$ must be at least 2. It follows that the Thurston semi-norm of $\omega = [F] \in H_2(M^{\varphi^N}_f;\R)$ is nonzero. 
\end{proof}


\section{Proof of Theorem \ref{thm_main}} \label{sec_proof}

\begin{thm_main}
		Let $\varphi$ be a self-diffeomorphism of $S_g$ that is not isotopic to the identity map. Then, there exists a manifold surface pair $(M,f)$ such that 
	$$\mathrm{rk} \bigoplus_{c_1(\mathfrak{s})\ne 0} \hfh(M,\mathfrak{s}) \ne \mathrm{rk} \bigoplus_{c_1(\mathfrak{s})\ne 0} \hfh(M_f^\varphi,\mathfrak{s})$$
	
\begin{proof}
	Let $G \vartriangleleft \mcg (S_g)$ be the set of $\hfh$-invisible mapping classes. We begin by showing that $G$ contains no pseudo-Anosov element of the Torelli group. Let $[\varphi]\in\mathcal{I}(S_g)$ be a pseudo-Anosov element of the Torelli group. Also let $(M,f)$ be a manifold surface pair  such that $M = S^1 \times S^2$ and for some $N\in \N$ the the mutant manifold $M_f^{\varphi^N}$ has a homology class with nonzero Thurston semi-norm. The existence of such a pair is guaranteed by Proposition \ref{prop_pA}.
	
	A simple computation shows that the Heegaard Floer homology of $M = S^1\times S^2$ is isomorphic to $\Z \oplus \Z$ and is supported entirely in the \spinc structure whose first Chern class is zero \citelist{\cite{Ozsvath2004b}*{\S 3}}.  Thus, the rank of the non-torsion summands of $\hfh(M)$ is zero
$$\mathrm{rk}\bigoplus_{c_1(\mathfrak{s})\neq 0}\hfh(M, \mathfrak{s}) = 0.$$

By construction $M_f^{\varphi^N}$ has a homology class with nonzero Thurston semi-norm. It follows from work of Hedden and Ni, that $\hfh(M_f^{\varphi^N})$ is nontrivial in at least one \spinc structure with nonzero first Chern class \citelist{\cite{Hedden2010}*{Thm. 2.2}}. In particular, the rank of the non-torsion summands is positive
$$\mathrm{rk}\bigoplus_{c_1(\mathfrak{s})\neq 0}\hfh(M_f^{\varphi^N}, \mathfrak{s}) > 0$$
and therefore
$$\mathrm{rk}\bigoplus_{c_1(\mathfrak{s})\neq 0}\hfh(M, \mathfrak{s}) \neq
\mathrm{rk}\bigoplus_{c_1(\mathfrak{s})\neq 0}\hfh(M_f^{\varphi^N}, \mathfrak{s}).$$

Thus, the mapping class $[\varphi^N] = [\varphi]^N$ is not $\hfh$-invisible. Because the $\hfh$-invisible mapping classes form a subgroup of $\mcg(S_g)$, we concluded that $[\varphi]$ is also not $\hfh$-invisible (Prop. \ref{prop_normal}). Therefore, no pseudo-Anosov element of the Torelli group is an element of $G$.

Furthermore, we showed in Propositions \ref{prop_normal} and  \ref{prop_genus-2}  respectively that $G$ is normal and does not contain the genus-2 hyperelliptic involution. Hence, $G$ is trivial by Proposition \ref{prop_MCGtheory}.
\end{proof}
\end{thm_main}

\section{Total rank detects mutation} \label{sec_totalrk}

\begin{thm_totalrk}
	Let $[\varphi]\in\mcg(S_g)$ be a mapping class that is isomorphic to neither the identity nor the genus-2 hyperelliptic involution. Then there exists a manifold surface pair $(M,f)$ such that
	$$\mathrm{rk} \ \hfh(M) \ne \mathrm{rk} \ \hfh(M_f^\varphi).$$

\begin{proof}
	Let $G$ be the set of mapping classes such that $[\varphi] \in G$ if and only if $\mathrm{rk} \ \hfh(M) = \mathrm{rk} \ \hfh(M_f^\varphi)$ for all manifold surface pairs $(M,f)$. Like the set of $\hfh$-invisible mapping classes $G$ is a normal subgroup of $\mcg(S_g)$. This follows from the proof of Proposition \ref{prop_normal} with the appropriate notation changes. Thus, it suffices to show that $G$ contains no pseudo-Anosov elements of the Torelli group (Proposition \ref{prop_MCGtheory}). 
		
	Let $[\varphi]\in\mathcal{I}(S_g)$ be a pseudo-Anosov element of the Torelli group. Also let $(M,f)$ be a manifold surface pair such that $M = S^1 \times S^2$ and for some $N\in \N$ the the mutant manifold $M_f^{\varphi^N}$ has a homology class with nonzero Thurston semi-norm. The existence of such a pair is guaranteed by Proposition \ref{prop_pA}.		
		
	Let $T$ be the result of 0-surgery on the trefoil. Hedden and Ni showed that $T$ and $M$ are the only closed, orientable, irreducible three-manifolds with nonzero first Betti number and $\mathrm{rk}\,\hfh = 2$ \citelist{\cite{Hedden2010}*{Thm. 1.1}}. In the proof of Proposition \ref{prop_pA}, we showed that the mutant $M_f^{\varphi^N}$ is closed, orientable and irreducible, and its first Betti number is nonzero. Thus, it is enough to show that $M_f^{\varphi^N}$ is not diffeomorphic to either $T$ or $M$.
		
	 A Mayer-Vietoris argument shows that $H_2(T;\R) \cong \R$. The Thurston semi-norm is constantly zero on $H_2(T;\R)$, because the trefoil is a genus-1 knot \citelist{\cite{Gabai1987a}*{Cor. 8.3}}. The homology group $H_2(M;\Z) = H_2(S^1\times S^2;\Z)$ is isomorphic to $\Z$  and is generated by the homology class of a sphere. Thus, the Thurston semi-norm of any homology class in $H_2(S^1\times S^2;\R)$ is zero. Therefore, the Thurston semi-norm differentiates $M_f^{\varphi^N}$ from both $T$ and $S^1\times S^2$. 	
\end{proof}
\end{thm_totalrk}

\section{Implications} \label{sec_implications}

There are two ways to interpret Theorems \ref{thm_main} and \ref{thm_totalrk} as statements about actions of mapping class groups of surfaces on categories. The first uses bordered Heegaard Floer homology and results in a statement about an action on a category of $\mathcal{A}_\infty$-modules. The second uses the definition of $\hfh$ and results in a statement about an action on a Fukaya category.

\subsection{Bordered Heegaard Floer homology} \label{sec_bordered}

In \citelist{\cite{Lipshitz2011a}} and \citelist{ \cite{Lipshitz2010a}}, Lipshitz, Ozsv\' ath and Thurston developed a variant of Heegaard Floer homology for three-manifolds with parametrized boundary called \emph{bordered Heegaard Floer homology}. These bordered invariants are related to $\hfh$ by pairing theorems\citelist{\cite{Lipshitz2010a}*{Thm. 1.3} \cite{Lipshitz2011a}*{Thm. 11}}. The pairing theorems provide a method for computing $\hfh(M)$ by cutting $M$ along separating surfaces and computing the bordered  Heegaard Floer homology of the resulting components. By applying this method to manifold surface pairs and their mutants, we can use Theorem \ref{thm_main} to infer information about the bordered Heegaard Floer homology of mapping cylinders of surface diffeomorphisms.

Let $\mcg_0(S_g)$ denote the \emph{strongly based mapping class group} of $S_g$ that is the isotopy classes of diffeomorphisms that fix a given disk in $S_g$. There is a canonical projection
$$p\colon \mcg_0(S_g) \to \mcg(S_g)$$
given by quotienting out by the copy of $\pi_1(S_g)$ that corresponds to pushing the disk around closed curves in $S_g$ as well as by the Dehn twist around the boundary of the disk. Following \citelist{\cite{Lipshitz2010a}*{\S 8}}, we assign to each strongly based mapping class $[\varphi] \in \mcg_0(S_g)$ the bimodule $\widehat{CFDA}(\varphi, 0)$ associated to its mapping cylinder. By considering Theorem \ref{thm_main} from the perspective of bordered Heegaard Floer homology, we get the following result about these bimodules.

\begin{cor} \label{cor_bordered}
	If $[\varphi] \in \mcg_0(S_g)$ is a strongly based mapping class such that $[\varphi]$ is not in the kernel of $p$, then the action of $[\varphi]$ on the category of $\mathcal{G}(Z)$-graded $\mathcal{A}(Z)$-modules given by tensoring with $\widehat{CFDA}(\varphi, 0)$ is not the trivial action.
	\begin{proof}
	Let $[\varphi] \in \mcg_0(S_g)$ such that $[\varphi]$ is not in the kernel of $p$. Also, let $(M,f)$ be a manifold surface pair such that the rank of the non-torsion summands of $\hfh(M)$ differs from that of $\hfh(M_f^\varphi)$. The existence of such a pair is guaranteed by Theorem \ref{thm_main}. Finally, let $M_1$ and $M_2$ be the connected components of $M\setminus f(S_g)$.
	
	 The Heegaard Floer homology of $M$ can be computed from the bordered invariants of $M_1$ and $M_2$ as follows
	$$\hfh(M) \cong H_*\left(\widehat{CFA}(M_1)\,\widetilde{\otimes}\,\widehat{CFD}(M_2)\right) $$
	where $\widetilde{\otimes}$ is the $\mathcal{A}_\infty$-tensor product.
	
	 Similarly, decomposing the mutant manifold $M_f^\varphi$ as the union $M_1 \cup C_\varphi \cup M_2$ where $C_\varphi$ is the mapping cylinder of $\varphi$ corresponds to the following module decomposition of $\hfh(M^\varphi_f)$.
		$$\hfh(M^\varphi_f) \cong H_*\left(\widehat{CFA}(M_1)\,\widetilde{\otimes}\,\widehat{CFDA}(\varphi, 0) \,\widetilde{\otimes}\, \widehat{CFD}(M_2)\right)$$
		
		Thus, the difference between $\hfh(M)$ and $\hfh(M_f^\varphi)$ must result from the effect of tensoring with $\widehat{CFDA}(\varphi,0)$. Therefore, the action of $[\varphi]$ on $\mathcal{A}(Z)$-modules given by tensoring with $\widehat{CFDA}(\varphi, 0)$ must not be the trivial action.
	\end{proof}
\end{cor}

A similar reformulation of Theorem \ref{thm_totalrk} gives the following result about the action of $\mcg_0(S_g)$ on the category of ungraded  $\mathcal{A}(Z)$-modules.

\begin{cor} \label{cor_bordered_totalrk}
	If $[\varphi] \in \mcg_0(S_g)$ is a strongly based mapping class such that $p([\varphi])$ is neither the identity nor the genus-2 hyperelliptic involution, then the action of $[\varphi]$ on the category of ungraded $\mathcal{A}(Z)$-modules given by tensoring with $\widehat{CFDA}(\varphi, 0)$ is not the trivial action.
\end{cor}

These results are closely related to work of Lipshitz, Ozsv\' ath and Thurston. In \citelist{\cite{Lipshitz2010}}, they showed that the action of a nontrivial strongly based mapping class $[\varphi]$ on the ungraded $\mathcal{A}(Z)$-modules given by tensoring with $\widehat{CFDA}(\varphi, \pm(g~-~1))$ is not the trivial action. 

\subsection{Fukaya categories} \label{sec_Fuk}

When viewed from another perspective, the work of Lipshitz, Ozsv\' ath and Thurston shows that the strongly based mapping class group $\mcg_0(S_g)$ acts freely on a version of the Fukaya category of $S_g$ with a disk removed as well as on a version of the Fukaya category of the $(2g-1)$-fold symmetric product $\mathrm{Sym}^{2g-1}(S_g - D)$ \cite{Auroux2010a}. 

Theorem \ref{thm_totalrk} is also related to mapping class group actions on Fukaya categories. In particular, the chain complex that underlies $\hfh$ of a three-manifold with a genus-$g$ Heegaard splitting corresponds to a morphism group in the Fukaya category of the $g$-fold symmetric product of $S_g$ with a point removed. Furthermore, the action of the based mapping class group of $S_g$ on the symmetric product $\mathrm{Sym}^g(S_g - z)$ induces a strict action on the Fukaya category $\mathrm{Fuk}(\mathrm{Sym}^g(S_g - z))$ \citelist{\cite{Seidel2008}*{\S 10b}}. 

\begin{cor} \label{cor_fukaya}
	If $[\varphi] \in \mcg(S_g - z)$ is a based mapping class such that the corresponding element of $\mcg(S_g)$ is neither the identity nor the genus-2 hyperelliptic involution, then the action of $[\varphi]$ on the Fukaya category $\mathrm{Fuk}(\mathrm{Sym}^g(S_g - z))$ is not the trivial action. In particular, the map induced by $\varphi$ on $\mathrm{Sym}^g(S_g - z)$ is not Hamiltonian isotopic to the identity.
	
	\begin{proof}
	Let $[\varphi] \in \mcg(S_g - z)$ be a based mapping class such that the corresponding element of $\mcg(S_g)$ is neither the identity nor the genus-2 hyperelliptic involution. Also, let $(M,f)$ be a manifold surface pair such that $f(S_g)$ is a Heegaard surface and 
	$$\mathrm{rk}\,~\hfh(M)\ne \mathrm{rk}\,~\hfh(M_f^\varphi).$$
	The existence of such a manifold is guaranteed by the fact that the proof of Theorem \ref{thm_totalrk} only uses manifold surface pairs where the embedded surface is a Heegaard surface. Finally, let $T_\alpha$ and $T_\beta$ be the corresponding Heegaard tori in $\mathrm{Sym}^g(S_g - z)$. 
	
	The action of $[\varphi]$ on $\mathrm{Fuk}(\mathrm{Sym}^g(S_g - z))$ sends $T_\beta$ to $T_{\varphi(\beta)}$, the Heegaard torus that results from translating the curves of $\beta$ by $\varphi$. Furthermore, $T_\alpha$ and $T_{\varphi(\beta)}$ are the Heegaard tori of a splitting of the mutant manifold $M_f^\varphi$. It then follows from the definitions that 
	$$\widehat{CF}(M) = \mathrm{Mor}(T_\alpha, T_\beta) \mbox{ and } \widehat{CF}(M_f^\varphi) = \mathrm{Mor}(T_\alpha, T_{\varphi(\beta)}).$$
	
	Because $\hfh(M)$ and $\hfh(M_f^\varphi)$ do not have the same rank, we concluded that their underlying chain complexes $\widehat{CF}(M)$ and $\widehat{CF}(M_f^\varphi)$ are not quasi-isomorphic. Thus, the morphism groups $\mathrm{Mor}(T_\alpha, T_\beta)$ and $\mathrm{Mor}(T_\alpha, T_{\varphi(\beta)})$ are not quasi-isomorphic. Therefore, $T_\beta$ is not isomorphic to $T_{\varphi(\beta)}$. 
	\end{proof}
\end{cor}

It should also be possible to reformulate Theorem \ref{thm_main} as a statement about an action of the based mapping class group of $S_g$ on a version of the Fukaya category of $\mathrm{Sym}^g(S_g - z)$. Such a reformulation would likely require working with grading data like that described in \cite{Sheridan2011}. We will return to this in a future paper.

\addcontentsline{toc}{section}{References}
\bibliography{DigLib}
\end{document}